\documentclass[12pt]{article}
\usepackage{amsmath, amssymb, amsfonts,bbm,bm}
\usepackage{graphicx}
\usepackage{enumerate}
\usepackage{natbib}
\usepackage{xurl} 
\usepackage{hyperref}
\usepackage{color}
\usepackage{setspace}
\newcommand{\blind}{1}

\usepackage[top=0.6in, bottom=0.7in,left=1in,right=1in]{geometry}
\usepackage{natbib}
 \bibpunct[, ]{(}{)}{,}{a}{}{,}%
\usepackage{rotating}
\usepackage{fancyvrb}
\usepackage{booktabs} 
\usepackage{ulem}
\usepackage{multirow} 
\usepackage{comment}

\usepackage{algorithm}
\usepackage{algorithmic}

\makeatletter
\newenvironment{breakablealgorithm}
{
\begin{center}
\refstepcounter{algorithm}
\hrule height.8pt depth0pt \kern2pt
\renewcommand{\caption}[2][\relax]{
 {\raggedright\textbf{\ALG@name~\thealgorithm} ##2\par}%
 \ifx\relax##1\relax 
   \addcontentsline{loa}{algorithm}{\protect\numberline{\thealgorithm}##2}%
 \else 
   \addcontentsline{loa}{algorithm}{\protect\numberline{\thealgorithm}##1}%
 \fi
 \kern2pt\hrule\kern2pt
}
}{
\kern2pt\hrule\relax
\end{center}
}
\makeatother

\bibpunct[, ]{(}{)}{,}{a}{}{,}%

\newtheorem{thm}{Theorem}

\newtheorem{assump}{Assumption}
\newtheorem{eg}{Example}

\def\norm1{{\psi}}
\def\faith{{\varsigma}}
\def\eigenvalue{{\mu}}

\newcommand{\RNum}[1]{\uppercase\expandafter{\romannumeral #1\relax}}

\newcommand{\Var}{{\mbox{Var}}}
\newcommand{\cov}{{\mbox{cov}}}

\allowdisplaybreaks

\begin{document}

\def\spacingset#1{\renewcommand{\baselinestretch}%
{#1}\small\normalsize} \spacingset{1}


\if1\blind
{
  \title{\bf Online Sparse Regression with Expanding Observables}
  \author{Ying Yang$^a$
    and 
    Fang Yao$^b$, 
    \thanks{
    Corresponding author: \href{fyao@math.pku.edu.cn}{fyao@math.pku.edu.cn} 
    }
    \vspace{.2cm}
    \\
    $^a$Center for Applied Mathematics, \\Shanghai Key Laboratory for Contemporary Applied Mathematics, \\Fudan University, Shanghai, China;\\
    $^b$School of Mathematics, Center for Statistical Science, \\Peking University, Beijing, China
    }
  \maketitle
} \fi

\if0\blind
{
  \bigskip
  \bigskip
  \bigskip
  \begin{center}
    {\LARGE\bf Online Sparse Regression with Expanding Observables}
\end{center}
  \medskip
} \fi

\bigskip
\begin{abstract}
Online high-dimensional regression has gained increasing attention in recent years, yet existing methods typically assume that all candidate features, including important ones, are observed from the outset of data collection. This assumption is often violated in real-world scenarios, where new variables become available gradually as data accumulate. To address this gap, we introduce a novel framework, Recurrent Adaptive Variable Selection (RAVAS), for online regression with expanding observability.  RAVAS employs a recurrent procedure that dynamically updates feature selection as both the sample size and the observable feature set grow. The algorithm is designed to be computationally efficient and memory-light, relying only on low-dimensional sufficient statistics that are updated online. A key advantage of the method lies in its ability to detect and incorporate important variables that emerge later, thereby mitigating the effect of early-stage missingness. We establish theoretical guarantees on model selection, estimation error, and feature coverage, and develop an adaptive online tuning strategy. Extensive simulations and real-world experiments verify the effectiveness of RAVAS for high-dimensional streaming data.
\end{abstract}

\noindent%
\vfill

\newpage
\spacingset{1.3} 

\section{Introduction}

High-dimensional statistical learning has been prevalent for decades due to its capacity to handle datasets with large numbers of features \citep{tibshirani1996regression,fan2001variable,zhang2010nearly,buhlmann2011statistics}. With advancements in data collection technologies, streaming data that are collected online with high speed have become increasingly popular. A key challenge that emerges at the intersection of high-dimensional statistics and real-time online learning is the dynamic expansion of data dimensionality during the collection process. For example, in air pollution monitoring, new environmental sensors are frequently deployed over time, expanding the range of observable variables. Likewise, in internet-based financial platforms such as peer-to-peer lending, additional borrower- or transaction-level attributes become available as platforms update their data pipelines, leading to a richer covariate set over time. Despite the growing importance of such dynamically evolving datasets, existing research has predominantly focused on the setting of fully observed feature space, where all candidate features, including important ones, are observed from the beginning of online data collection. In this work, we investigate online estimation in high-dimensional regression under expanding observables, where both the sample size and the number of observed features increase over time.

\subsection{Related work}
\label{sec:literature}

A significant body of research has focused on online convex optimization with regularization. \cite{foster2016online} proposed an algorithm that inputs a limited number of features in each round to minimize the regret function. \cite{langford2009sparse} introduced sparse online learning via truncated gradient. \cite{sun2021novel} presented a framework for online learning based on running averages. Additionally, \cite{duchi2009efficient} proposed the forward-backward splitting method, establishing a framework for online proximal gradient. Other related works include \cite{zhao2020unified,vijayakumar2005incremental,ma2017stabilized,ning2020control}. Beyond convex optimization, studies have also explored online feature selection \citep{ma2020variable,zhou2017online,wu2021latent}, discrete optimization \citep{xu2013adaptive} and  matrix factorization \citep{mairal2010online}. These studies emphasize the optimization performance of the algorithms, specifically, the convergence of online algorithms compared to classical ones given the same regularization parameter. However, they do not address the selection of the regularization parameter or the issue of statistical convergence.

To address this gap, statisticians have proposed several methods. \cite{han2021online} studied an online version of the debiased lasso and introduced the ``rolling-original-recalibration'' method to select the tuning parameter, using historical statistics as the training set and newly observed data as the testing set. \cite{yang2023online} solved an optimization problem that incorporates a linearized loss function with a pre-specified regularization parameter. \cite{fan2018statistical} suggested a two-step approach: first, a burn-in step to obtain an initial estimate $\widehat\beta_0$, followed by a refinement step using a variant of truncated stochastic gradient descent (SGD) on the support of $\widehat\beta_0$. Beyond regression, online algorithms have gained prominence in various high-dimensional statistical applications, including dynamic pricing \citep{fan2024policy, wang2024pricing}, sequential decision-making \citep{bastani2020online, wang2024online} and high dimensional reinforcement learning \citep{hao2021online}.

\subsection{Expanding observable regime and challenges}

\begin{figure}[htp]
\centering
\includegraphics[width=1\textwidth]{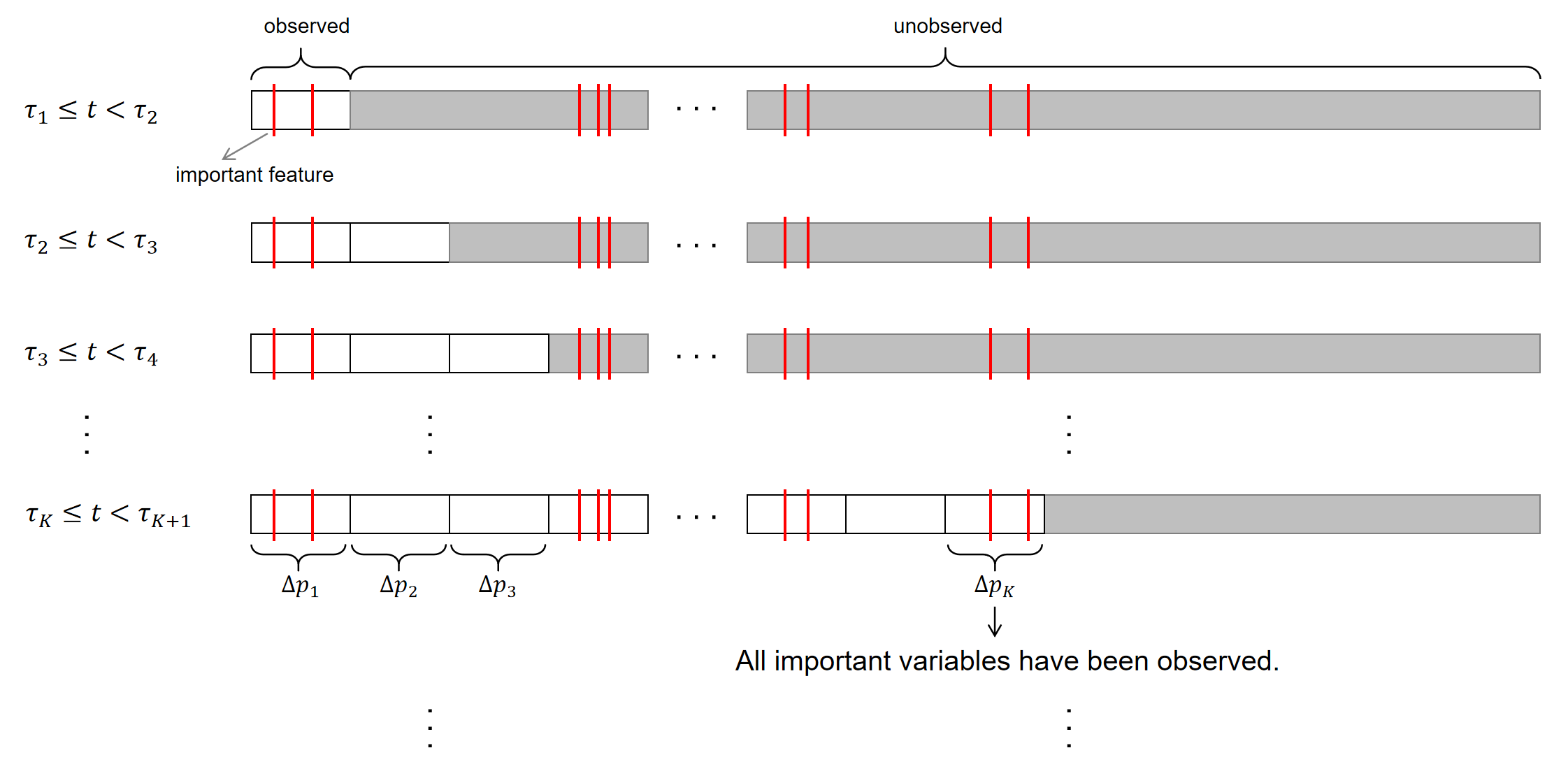}
\vspace{-15pt}
\caption{\label{fig:illustration} An illustration of the data observation process. The white segments represent the observed features and the gray ones represent the unobserved variables. The red lines correspond to the underlying important features. After time $\tau_K$, all key variables become observable.}
\end{figure}

In this work, we focus on the scenarios where both the sample size  and the observed feature dimension  increase over time. As illustrated in Figure \ref{fig:illustration}, 
at time $t$, we observe the data block  $\mathcal{D}_t=\{X_t,y_t\}$ containing $n$ observations, where $X_t=(X_{t1}^\top,\ldots,X_{tn}^\top)^\top\in\mathbb{R}^{n\times p_t}$ is the observed design matrix and $y_t=(y_{t1},\ldots,y_{tn})^\top\in\mathbb{R}^{n}$ is the response vector. We assume that the underlying data-generating process is  
\begin{equation*}
y_t = X_t \alpha_*^{(t)} + U_t\gamma_*^{(t)} + \varepsilon_t,
\end{equation*}
where  $U_t\in\mathbb{R}^{n\times q_t}$ is the unobserved variables, $\alpha_*^{(t)}\in\mathbb{R}^{p_t}$ and $\gamma_*^{(t)}\in\mathbb{R}^{q_t}$ are the regression coefficients, and $\varepsilon_t$ is the $n$-dim noise vector. We assume that $\alpha_*^{(t)}$ is sparse, satisfying $|J(\alpha_*^{(t)})|=s_t\ll p_t$. Denote $Z_{ti}=(X_{ti}^\top, U_{ti}^\top)^\top\in\mathbb{R}^{p_Z}$, which is the full feature space with dimension $p_Z$. 
As data collection progresses, more features become observable, meaning that $p_t$ increases while 
$q_t$ decreases., and the total number of potential features is $p_Z$.
We focus on the scenario in which the underlying model is characterized by the existence of a vector $\beta_*$ such that for any $t$, $\beta_*=(\alpha_*^{(t),\top},\gamma_*^{(t),\top})^\top$,  and
\begin{equation}
\label{eq:full model}
y_t = Z_t \beta_* + \varepsilon_t.
\end{equation}
{Such data observation process is significantly distinct with the settings explored in existing works, 
where the full feature space is assumed to be observable from the outset.}

As the observed model dimension $p_t$ grows, the memory and computational costs of storing and updating high-dimensional covariance statistics increase rapidly, eventually exceeding the limits of computational resources. More importantly, online algorithms aim to operate within the time and memory budgets; retaining all features at all times is therefore infeasible. Hence, dynamically selecting and maintaining only a small subset of influential variables becomes essential to ensure statistical tractability and computational efficiency over time.

Conducting variable selection in the increasing-observable case presents four primary challenges.
Firstly, features are increasingly observed, which may leave certain important features unobserved at the initial stage. This raises the question of how to implement variable selection while ensuring theoretical guarantee that the estimated model covers the observed important features with high probability. 
Secondly, as the observed dimension dynamically increases, it is necessary to ensure that the variable selection process efficiently reduces the model dimension to keep computational costs acceptable. 
Thirdly, the accumulation of estimation error over multiple variable selection operations complicates the analysis, requiring a systematic approach to quantify this error.
Finally, selecting the appropriate tuning parameter remains a significant challenge that can impact the overall effectiveness of the variable selection process.

\subsection{Our contributions}

To address the above challenges, we establish a novel framework for online high-dimensional regression with expanding observables. 
For the first challenge, we introduce the regression faithfulness and seperation structure assumptions to ensure the feasibility of high-dimensional regression in the presence of unobserved variables. When certain important variables are unavailable, the conditional correlation analysis between the response and any feature can be biased, potentially leading to the erroneous exclusion of important features. Inspired by the faithfulness assumption in causal discovery, we propose the concept of regression faithfulness to mitigate this issue.

With respect to the second challenge, we propose the Recurrent Adaptive Variable Selection (RAVAS) method. RAVAS operates through a recurrent adaptive procedure that updates the active model as new features become observable. Each update cycle begins with an initial estimator, followed by a soft selection stage in which a lasso estimator with a deliberately relaxed penalty is applied. This enables broad yet probabilistically reliable identification of potentially important variables under limited effective sample size. 
When sufficient data have accumulated, the procedure may further transition to a hard selection stage, yielding a more compact model and improved estimation accuracy. As new features continue to emerge, RAVAS is re-applied to the expanded design, ensuring that the model dimension remains tractable and that important predictors are dynamically identified throughout the process.

To address the third challenge, we develop a theoretical framework that characterizes the behavior of the dynamically selected model throughout the online procedure. We show that RAVAS achieves reliable dimension reduction at each cycle and, with high probability, retains all truly relevant variables. Additionally, we quantify the estimation error, showing that the method achieves low error with high probability. 
Finally, to enhance practical implementation, we introduce a data-driven online tuning strategy. By extending existing techniques to the online setting, we consistently estimate the residual variance, while also adapting to the complexities introduced by unobserved variables. An online cross-validation procedure is introduced to efficiently estimate the tuning parameter, facilitating dynamic updates of sufficient statistics as data accumulate.

The remainder of the paper is organized as follows. Section \ref{sec:model} presents the RAVAS method, detailing its algorithmic implementation and the proposed online cross-validation procedure.
\ref{sec:theory} first introduces key assumptions, including regression faithfulness and separation structure assumptions, and then establishes theoretical guarantees for the method. Technical proofs of main theorems are deferred to the Supplementary Material. 
Sections \ref{sec:simu} and \ref{sec:realdata} validate our findings through numerical simulations and real data applications, respectively.

\textbf{\textit{Notations}.} For an arbitary vector $a\in\mathbb{R}^m$, define the norm $\|a\|_q=\{\sum_{i=1}^m |a_i|^q\}^{1/q}$ where $q\in\mathbb{N}$.  
For a vector $\beta\in\mathbb{R}^p$, let $J(\beta)=\{j:\beta_j\neq0\}$  represent the index set of nonzero elements of $\beta$. For an index set $J\subset\{1,2,\ldots,p_t\}$, $|J|$ is the cardinal of $J$, $J^c$ is the complementary set $J$ in $\{1,2\ldots,p_t\}$ and $\beta_J$ is the subvector of $\beta$ with elements indexed in $J$. For two positive quantities $a$ and $b$, the notation $a\gtrsim b$ means $b/a=O(1)$.

\section{Recurrent Adaptive Variable Selection}
\label{sec:model}

In model \eqref{eq:full model}, we assume that $\{Z_{ti}\}_{t,i}$ are independently and identically distributed (i.i.d.) random vectors with mean zero and covariance  $\Sigma_Z=E Z_{\iota}^\top Z_{\iota}$.
To describe changes of the observed data, we introduce  several key quantities. The change points in the observed dimension are defined as $\mathcal{T}^p=\{\iota:\ p_\iota-p_{\iota-1}>0\}=\{\tau_{k}:k=1,\ldots,K_m\}$.  Let $\Delta p_k=p_{\tau_{k}}-p_{\tau_k-1}$ represent the increase in observed dimension at time $\tau_{k}\in \mathcal{T}^p$ with $p_0=0$ and $\Delta p=\max_k \Delta p_k$ be the maximum dimension increment. Obviously, $\tau_1=1$ and $\Delta p_1=p_1$. We assume that there exists a $K\le K_m$ such that for $t<\tau_K$,  certain important variables are unobserved, which corresponds to $J(\gamma_*^{(t)})\neq \emptyset$; conversely, for $t\ge\tau_K$, all important features are present, indicating $\gamma_*^{(t)}=0$, as illustrated in Figure \ref{fig:illustration}. However, the value of $K$ is unknown. 

In the expanding-observable scenario, the key to keep the computational burden affordable throughout is to dynamically select important features. One related method is variable screening, which narrows the selection range by screening out unimportant features. Existing variable screening methods can be divided into two categories: those based on marginal correlation and those based on conditional correlation. Marginal-correlation methods require strong assumptions about the correlations between each variable and the outcome, making them restrictive in practice \citep{fan2008sure,fan2010generalized}. Conditional-correlation methods, on the other hand, calculate dependencies using all covariates at each step, resulting in high computational costs \citep{wang2009forward,buhlmann2010variable,cho2012high,wang2016high}. Therefore, a new dimension reduction technique is needed that balances effectiveness with computational efficiency in the expanding observables setting. 
Motivated by the observation that the probability of the lasso-estimated model covering the true model increases as the penalty parameter decreases, we aim to develop a lasso-based variable selection method with a data-driven tuning parameter that adapts seamlessly to both the incomplete- and complete-observable contexts though $\tau_K$ is unknown.

\subsection{Recurrent adaptive variable selection (RAVAS)}

With slight abuse of notation, for any $t$, let $\tau_k=\max\{\tau\in\mathcal{T}^p:\tau \le t\}$ represent the most recent change point before time $t$. The number of effective batches of observations that have the most complete observability up to $t$ is then given by
$$\zeta_t=t-\tau_k+1.$$
Earlier models deviate further from the true model, and incorporating observations with a higher proportion of unobserved variables offers no advantage. Thus, the efficient strategy is to estimate the regression using the $\zeta_tn$ most recent effective samples.
Historical data, however, remains essential for guiding variable selection: it enables us to identify and retain potentially important features before estimation, thereby reducing the model dimension and ensuring that subsequent computation remains feasible. When $\zeta_t$ is small and $\Delta p_{\tau_k}$ is large,  this scenario remains a ``small $n$, large $p$'' problem.
The lasso problem based on data $\{(X_i,y_i)\}_{i=\tau_k}^t$ ($t<\tau_{k+1}$) is 
\begin{equation*}
\arg\min_\alpha \frac{1}{\zeta_tn}\sum_{i=\tau_k}^t\|y_i-X_i\alpha\|_2^2+\lambda\|\alpha\|_1.
\end{equation*}
where $\lambda$ is the regularization parameter. This is equivalent to solve 
\begin{equation}\label{eq:naive taget}
\arg\min_\alpha \left(-2C_{Xy}^{(t)}\alpha+\alpha^\top C_{X}^{(t)} \alpha\right)+\lambda\|\alpha\|_1,
\end{equation}
where $\ C_{Xy}^{(t)}=(\zeta_tn)^{-1}\sum_{i=\tau_k}^tX_i^\top y_i\ \text{and}\ C_{X}^{(t)}=(\zeta_tn)^{-1}\sum_{i=\tau_k}^tX_i^\top X_i$.
Note that $C_{y}^{(t)},C_{Xy}^{(t)},C_{X}^{(t)}$ can be updated in an online fashion:
{\begin{align*}
&C_{Xy}^{(t)}=\frac{\zeta_t-1}{\zeta_t}\cdot C_{Xy}^{(t-1)}+\frac{1}{\zeta_tn}X_{t}^\top y_{t},\ C_{X}^{(t)}=\frac{\zeta_t-1}{\zeta_t}\cdot C_{X}^{(t-1)}+\frac{1}{\zeta_tn}X_{t}^\top X_{t}.
\end{align*}}
The naive method is to treat $\{C_{Xy}^{(t)},C_{X}^{(t)}\}$ as sufficient statistics. However, such methods are not computationally  efficient as storing the statistics $\{C_{Xy}^{(t)},C_{X}^{(t)}\}$ and solving \eqref{eq:naive taget} require the computational complexity $O(p_t^2)$, making it inefficient for high-dimensional settings.

To overcome this difficulty, we propose the RAVAS method, which dynamically selects important
features with high probability through a structured three-stage procedure.  
Let $\Delta p$ denote the maximal dimension increment. We set the relative timestamps for the
soft and hard selection stages as
\[
\iota^0 = (\log \Delta p)^\kappa / n, 
\qquad 
\iota^* = c_h \Delta p / n,
\]
where $\kappa>1$ and $c_h$ is a user-specified constant. A larger value of $c_h$ (e.g.\ $c_h\ge 1$)
yields a more conservative strategy that increases the likelihood of retaining important variables,
while a smaller value accelerates dimensionality reduction and reduces computational cost.
Numerical experiments suggest $c_h \in [0.3,1]$ as a practical range.

Let $d_t$ denote the dimension of the selected feature set at time $t$ (distinct from the ambient
dimension $p_t$). For $\tau_k \le t < \tau_{k+1}$, the estimation procedure typically proceeds
through the following stages.
\begin{itemize}
\item[1.] \textit{\textbf{Initial stage}} (Warm-up stage). The initial stage serves to accumulate a sufficient number of effective samples while simultaneously providing a real-time preliminary estimate when needed in practice.
For time $\zeta_t \le \iota^0$, no variable selection is conducted, and the estimated index set is
$$\widehat{J}_{t}=\widehat{J}_{\tau_k}=J(\widehat\alpha_{\tau_k-1})\cup\{p_{\tau_k-1}+1,\ldots,p_{\tau_k}\},
$$
which includes indices of important variables at time $\tau_k-1$ along with indices of newly emerged features at time $\tau_k$. The estimated model dimension is given by $d_t=d_{\tau_k}=|\widehat{J}_{\tau_k}|$. 
We use the classical methods to solve
\begin{equation}
\label{eq:original target increase}
\frac{1}{\zeta_tn}\sum_{\iota=\tau_k}^t(y_\iota-X_{\iota,\widehat{J}_{\tau_k}}\alpha)^2+\lambda_t\|\alpha\|_1
\end{equation} 
with the penalty parameter $\lambda_t=\lambda_t^*=O\left(\sqrt{\log d_{\tau_k}/(\zeta_tn)}\right)$
to obtain the estimator $\widehat{\alpha}_{t}$. When $\zeta_t=\iota^0$, we compute the initial statistics of this cycle
\begin{equation}\label{eq:initial stat increase}
C_{Xy}^{(t)}(\widehat{J}_t)=\frac{1}{\zeta_tn}\sum_{\iota=1}^tX_{\iota,\widehat{J}_{t}}^\top y_{\iota},\quad
C_{X}^{(t)}(\widehat{J}_t)=\frac{1}{\zeta_tn}\sum_{\iota=1}^tX_{\iota,\widehat{J}_{t}}^\top X_{\iota,\widehat{J}_{t}}.
\end{equation}
and remove the stored observations. 

\item[2.]  \textit{\textbf{Soft selection} ($\iota^0<\zeta_t\le\min\{\tau_{k+1}-\tau_k,  \iota^*\}$)}. Since the probability of covering the true model with the lasso method increases as the tuning parameter $\lambda$ decreases, we perform lasso regression with a penalty parameter $\lambda_t^0$ that is smaller than the optimal value to estimate the support of important variables. Specifically, for $\iota^0<\zeta_t\le\iota^*$, with the prior constraint $\widehat{\alpha}_{t, \widehat{J}_{t-1}^c}=0$, we solve the $d_{t-1}=|\widehat{J}_{t-1}|$-dimension lasso problem
\begin{equation}\label{eq:online lasso raoa}
\widehat{\alpha}_{t, \widehat{J}_{t-1}}=\arg\min_\alpha \left\{-2 {C}_{Xy}^{(t)}(\widehat{J}_{t-1}) \alpha+\alpha^\top {C}_{X}^{(t)}(\widehat{J}_{t-1})\alpha+\lambda_t^0\|\alpha\|_1\right\},
\end{equation}
where $\lambda_t^0$ is set smaller than the optimal tuning parameter, {which is explicitly defined in  \eqref{eq:lam0 raoa} in the following subsection. The optimization problem in \eqref{eq:online lasso raoa} is solved using the proximal gradient algorithm, detailed in Algorithm 3 in the Supplementary Material.} Then we update $\widehat{J}_{t}=J(\widehat{\alpha}_{t})$ and the statistics, 
\begin{align}\label{eq:update stat raoa}
&C_{Xy}^{(t)}(\widehat{J}_{t})=\frac{\zeta_t-1}{\zeta_t}C_{Xy}^{(t-1)}(\widehat{J}_{t})+\frac{1}{\zeta_tn}X_{t,\widehat{J}_{t}}^\top y_{t,i},\quad C_{X}^{(t)}(\widehat{J}_{t})=\frac{\zeta_t-1}{\zeta_t}C_{X}^{(t-1)}(\widehat{J}_{t})+\frac{1}{\zeta_tn}X_{t,\widehat{J}_{t}}^\top X_{t,\widehat{J}_{t}}.
\end{align}
We demonstrate in Section \ref{sec:theory} that, under mild conditions, the selected features can capture the observed important variables with high probability. 

\item[3.]  \textit{\textbf{Hard selection} ($\zeta_t > \iota^*$)}. As the number of selected important features $d_t$ reduces and the effective sample size $\zeta_tn$ grows, the problem gradually approaches an ordinary regression scenario. In this case, the procedure may advance to the hard selection stage, where variable selection is performed using a least squares–based criterion to further refine the model dimension.
If the effective sample size is insufficient before the next batch of variables appears, the algorithm simply proceeds to the next cycle; nonetheless, the soft-selection stage alone guarantees adequate dimensionality control for subsequent updates.
When $\zeta_t > \iota^*$, we calculate the OLS estimator $\widehat{\alpha}_t^{ols}=\{C_{X}^{(t)}(\widehat{J}_{t-1})\}^{-1}C_{Xy}^{(t)}(\widehat{J}_{t-1})$ and then estimate $\alpha_*^{(t)}$ by the following hard selection,
\begin{equation}\label{eq:hard}
\widehat{\alpha}_{t}=\text{thre}\left(\widehat{\alpha}_{t}^{ols}\right),\text{ where }\text{thre}\left(\widehat{\alpha}_{t}^{ols}\right)_j=
\left\{
\begin{aligned}
&0, \ \text{if}\ |\widehat{\alpha}_{t,j}^{ols}|<b_j;\\ 
&\widehat{\alpha}_{t,j}^{ols},\ \text{else}. 
\end{aligned}\right.
, 
\end{equation}
where  {$b_j=c_b(\iota^0n)^{-1/2}\widehat{\sigma}_e$} with a constant $c_b>0$ and the noise level estimator $\widehat{\sigma}_e$ detailed in the next subsection.  Finally, the index set of estimated important variables is updated as $\widehat{J}_t=J(\widehat\alpha_t)$. This transition is motivated by the advantages of streaming data, where the sample size increases over time, and the preceding selection stage helps to gradually lower the model dimension. Eventually, the Gram matrix will become positive definite.
\end{itemize}


We illustrate the procedure of RAVAS method in Figure \ref{fig:procedure} and outline the detailed implementation steps in Algorithm \ref{alg:RAVAS} below. 
Note that $\lambda_t^0$ is intentionally set at a relaxed level to ensure the inclusion of important variables with high confidence. To achieve more efficient estimation, we propose using an alternative penalty, $\lambda_t^*$, which we set at the optimal order at each step and defer the discussion of selecting $\lambda_t^0$ and $\lambda_t^*$ to the next subsection. 

\textit{Remark.}
Although RAVAS is developed for settings with expanding observables, the procedure also applies to the degenerate fully observable case in which all covariates are available from the beginning \citep{fan2018statistical,han2021online,yang2023online,wang2024pricing}. In this scenario, the algorithm reduces to a single-cycle implementation with the same soft–hard selection mechanism, and the theoretical guarantees continue to hold. We provide further discussion and numerical results for this special case in Section 3 of the Supplementary Materials.

\begin{figure}[H]
\centering
\includegraphics[width=1\textwidth]{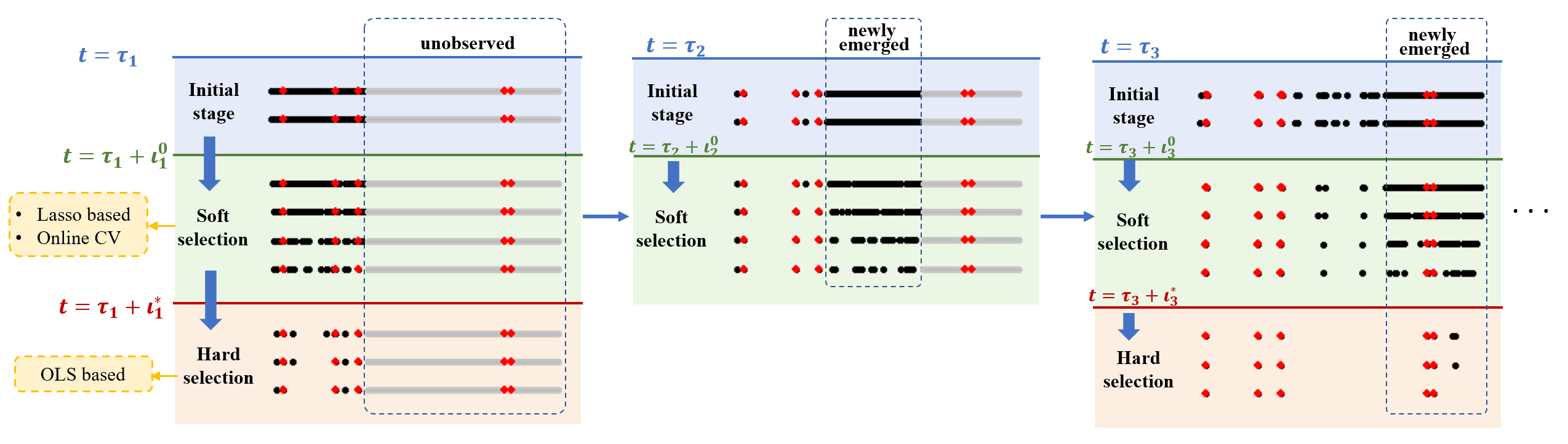}
\vspace{-0.5in}
\caption{\label{fig:procedure} An illustrate for the procedure of RAVAS method. }
\end{figure}

\begin{onehalfspace}
\begin{breakablealgorithm}
  \renewcommand{\algorithmicrequire}{\textbf{Input:}}
  \renewcommand{\algorithmicensure}{\textbf{Output:}}
  \caption{Recurrent Accelerating Online Algorithm}\label{alg:RAVAS}
  \begin{algorithmic}[1]
  \STATE \textbf{Initialize:} $\mathcal{O}^{(t)}\leftarrow\emptyset$; $p_0 \leftarrow 0$; $\widehat{J}_0 \leftarrow \emptyset$.  
  \STATE \textbf{For} $t=1,2,\ldots, T_{max}$:
  \STATE \quad \textbf{input}: $\{X_t, y_t\}$ and $p_t = \text{dim}(X_t)$
  \STATE \quad \textbf{if} $p_t>p_{t-1}$:
  \STATE \quad\quad $\zeta_t \leftarrow 1$;
  \STATE \quad\quad $\widehat{J}_{t-1}\leftarrow \widehat{J}_{t-1}\cup\{p_{t-1}+1, \ldots, p_t\}$, $d_{t-1} \leftarrow |\widehat{J}_{t-1}|$; 
  \STATE \quad \textbf{else}:
  \STATE \quad\quad $\zeta_t \leftarrow \zeta_t + 1$;
  \STATE \quad \underline{\textbf{if} $\zeta_t \le \iota^0$ (initial stage)}:
  \STATE \quad \quad $\mathcal{O}^{(t)}\leftarrow\mathcal{O}^{(t-1)}\cup \{X_{t,\widehat{J}_{t-1}}, y_t\}$;
  \STATE \quad \quad solve $\widehat{\alpha}_{t}$ based on $\mathcal{O}^{(t)}$ by \eqref{eq:original target increase} with $\lambda_t^*=2\sqrt{\log d_{t-1} /(\zeta_tn)}$;
  \STATE \quad \quad \textbf{if} $\zeta_t=\iota^0$:
  \STATE \quad \quad \quad solve $\widehat{\alpha}_{t}$ according to \eqref{eq:original target increase} with $\lambda_t^0$ defined in \eqref{eq:lam0 raoa};
  \STATE \quad \quad \quad initialize statistics $C_y^{(t)},C_{Xy}^{(t)}(\widehat{J}_t),C_{X}^{(t)}(\widehat{J}_t)$ by \eqref{eq:initial stat increase};
  \STATE \quad \quad \quad $\mathcal{O}^{(t)}\leftarrow\emptyset$;
  \STATE \quad \underline{\textbf{if} $\iota^0\le \zeta_t\le \iota^*$ (soft selection stage)}: 
  \STATE \quad \quad update $\lambda_t^0$ by \eqref{eq:lam0 raoa} and solve $\widehat{\alpha}_t$ according to \eqref{eq:online lasso raoa} with $\lambda_t^0$;
  \STATE \quad \quad update $\lambda_t^*$ by \eqref{eq:lam opt raoa} and solve $\widehat{\alpha}_t^*$ according to \eqref{eq:online lasso raoa} with $\lambda_t^*$;
  \STATE \quad \quad update $\widehat{J}_t\leftarrow J(\widehat{\alpha}_{t})$ and the statistics by \eqref{eq:update stat raoa};
  \STATE \quad \underline{\textbf{if} $\zeta_t > \iota^*$ (hard selection stage)}: 
  \STATE \quad\quad compute $\widehat{\alpha}_t^{ols}$, update $\widehat{\alpha}_t^*$, $\widehat{\alpha}_t\leftarrow \text{thre}(\widehat{\alpha}_t^{ols})$ by \eqref{eq:hard} and $\widehat{J}_t\leftarrow J(\widehat{\alpha}_{t})$;
  \end{algorithmic}  
\end{breakablealgorithm}
\end{onehalfspace}

\subsection{Online cross-validation and noise level estimation}
We now address the selection of $\lambda_t^0$. 
The residual of regressing $Y_t$ on $X_t$ is
\begin{equation}  
\label{eq:residual}
e_{ti}=\varepsilon_{ti}+U_{ti}^\top\gamma_*^{(t)}-X_{ti,\widehat{J}_{t-1}}^\top\Sigma_{X_{\widehat{J}_{t-1}}}^{-1}\Sigma_{X_{\widehat{J}_{t-1}}U}\gamma_*^{(t)}
\end{equation}  
and let $\sigma_e^2=\Var(e_{ti})$ be its variance. Let $C_\lambda$ be a positive constant that depends on the sub-exponential norm of the product of covariate and residual. Then the penalty parameter is selected as follows,
\begin{equation}
\label{eq:lam0 raoa}
\lambda_t^0= C_\lambda\sigma_e\sqrt{\frac{(\log d_{t-1})^{\delta_t}}{\zeta_tn}}\text{ with }\delta_t=\min\left\{1, \frac{\log(\zeta_tn)}{\log d_{t-1}}\right\}.
\end{equation}
Such $\lambda_t^0$ helps remove unimportant variables. Concurrently, to obtain the optimal estimator $\widehat\beta_t^*$, we use the following penalty parameter 
\begin{equation}
\label{eq:lam opt raoa}
\lambda_t^*=C_\lambda\sigma_e\sqrt{\frac{\log d_{t-1}}{\zeta_tn}}.
\end{equation}
Our selection of $\lambda_t^0$ ensures $\lambda_t^0 \ll 1$ even when $\log d_{t-1}\gg \zeta_tn$. 
When $\zeta_tn> d_{t-1}$, this approach reduces to the classical lasso setting. 

When there are unobserved important variables, the residual $e_{ti}$ in \eqref{eq:residual}  becomes dependent on the covariate $X_{ti,\widehat{J}_{t-1}}$, making it challenging to explicitly express $C_\lambda$. To address this, we propose an online version of $L$-fold cross validation to select $C_\lambda$. Let $C_y^{(t)}=(\zeta_tn)^{-1}\sum_{\iota=\tau_k}^ty_t^\top y_t$ be the sample variance of $y$, which can be updated onlinely similar to formula \eqref{eq:update stat raoa}.
Instead of calculating and storing $C_y^{(t)}$, $C_X^{(t)}(\widehat{J}_{t})$ and $C_{Xy}^{(t)}(\widehat{J}_{t})$, we need store $L$ sets of statistics, denoted by $C_y^{(t,\ell)}$, $C_X^{(t,\ell)}(\widehat{J}_{t})$ and $C_{Xy}^{(t,\ell)}(\widehat{J}_{t})$ for $\ell=1,2,\ldots,L$, where it holds that $C_y^{(t)}=\sum_{\ell=1}^{L}C_y^{(t,\ell)}$, $C_X^{(t)}(\widehat{J}_{t})=\sum_{\ell=1}^{L}C_X^{(t,\ell)}(\widehat{J}_{t})$ and $C_{Xy}^{(t)}(\widehat{J}_{t})=\sum_{\ell=1}^{L}C_{Xy}^{(t,\ell)}(\widehat{J}_{t})$. At each time $t$, we split the data block $\mathcal{D}_t$ into $L$ folds, denoted as $\mathcal{D}_t^{[1]}, \mathcal{D}_t^{[2]}, \ldots, \mathcal{D}_t^{[L]}$. Then we update  the $\ell$-th set of statistics based on the data set $\mathcal{D}_t^{[\ell]}$. 
For a value $C_{\lambda,j}$ from the given set for $C_\lambda$, we aim to use $\cup_{s=\tau_k}^t (\mathcal{D}_s/\mathcal{D}_s^{[\ell]})$ as the training set and $\cup_{s=\tau_k}^t \mathcal{D}_s^{[\ell]}$ as the test set. Thanks to the parametric form of the model, this can be achieved using the statistics $\{C_y^{(t,\ell)},C_X^{(t,\ell)}(\widehat{J}_{t}),C_{Xy}^{(t,\ell)}(\widehat{J}_{t})\}_{\ell=1}^{L}$. Specifically, let 
$${C}_{y}^{(t,-\ell)}=C_y^{(t)}-C_y^{(t,\ell)},\ {C}_{X}^{(t,-\ell)}(\widehat{J}_{t})=C_X^{(t)}(\widehat{J}_{t})-C_X^{(t,\ell)}(\widehat{J}_{t}),\ \text{and}\ {C}_{Xy}^{(t,-\ell)}(\widehat{J}_{t})=C_{Xy}^{(t)}(\widehat{J}_{t})-C_{Xy}^{(t,\ell)}(\widehat{J}_{t}).$$
and $\lambda_j=C_{\lambda,j}\sigma_e\sqrt{L(\log d_{t-1})^{\delta_t}/(L-1)\zeta_tn}$, we solve 
$$\widehat{\alpha}^{[\ell]}(\lambda_j)=\arg\min_\alpha \left\{{C}_{y}^{(t,-\ell)}-2 {C}_{Xy}^{(t,-\ell)}(\widehat{J}_{t}) \alpha+\alpha^\top {C}_{X}^{(t,-\ell)}(\widehat{J}_{t})\alpha\right\}+\lambda_j\|\alpha\|_1.
$$
The corresponding test error is  
$$err_j^{[\ell]}={C}_{y}^{(t,\ell)}-2 {C}_{Xy}^{(t,\ell)}(\widehat{J}_{t}) \widehat{\alpha}^{[\ell]}(\lambda_j)+\widehat{\alpha}^{[\ell]}(\lambda_j)^\top {C}_{X}^{(t,\ell)}(\widehat{J}_{t})\widehat{\alpha}^{[\ell]}(\lambda_j).
$$
We select $C_\lambda$ by setting $C_\lambda=C_{\lambda,j^*}$ where
$j^*=\arg\min_j \sum_{\ell=1}^{L} err_j^{[\ell]}$. 


Next, we address the online estimation of the noise level $\sigma_e$.  During the soft selection stage, similar to the debiased lasso method proposed by \cite{zhang2017simultaneous}, we use the scaled lasso approach proposed by \cite{sun2012scaled}. The noise level $\sigma_e$ is estimated via the following iterative algorithm: 
\begin{align*}
& \widehat{\sigma}_{e}^{(t)} \leftarrow (\zeta_tn)^{1/2} \left\{ C_{y}^{(t)}-2C_{Xy}^{(t)}(\widehat{J}_{t})\widehat{\alpha}_{\text{old}}+\widehat{\alpha}_{\text{old}}^\top C_{X}^{(t)}(\widehat{J}_{t}) \widehat{\alpha}_{\text{old}}\right\}^{1/2} , \\
& \lambda \leftarrow \widehat{\sigma}_{e}^{(t)}  \lambda_0, \ \widehat{\alpha}_{\text {new }}\leftarrow \min_\alpha L_\lambda\left({\alpha}\right), \\
& \widehat{\alpha}_{\text{old}} \leftarrow \widehat{\alpha}_{\text {new }}, 
\end{align*}
where $\lambda_0$ is a constant and
$$L_\lambda(\alpha)=\frac{C_{y}^{(t)}-2C_{Xy}^{(t)}(\widehat{J}_{t})\alpha+\alpha^\top C_{X}^{(t)}(\widehat{J}_{t}) \alpha}{\zeta_tn}+\lambda\|\alpha\|_1.
$$
This is equivalent to solving
\begin{equation*}
\arg\min_{\alpha,\sigma_\varepsilon}\left\{\frac{C_{y}^{(t)}-2C_{Xy}^{(t)}(\widehat{J}_{t})\alpha+\alpha^\top C_{X}^{(t)} (\widehat{J}_{t})\alpha}{\zeta_tn\sigma_\varepsilon}+\sigma_\varepsilon+\lambda\|\alpha\|_1\right\}.
\end{equation*}
As discussed in \cite{sun2012scaled} and \cite{sun2013sparse}, the universal penalty level $\lambda=\{2\log d_{t-1}/\zeta_tn\}^{1/2}$  generally performs well numerically. Additionally, \cite{sun2013sparse} propose selecting $\lambda=\sqrt{2}L_t(k/d_{t-1})$ with $k=L_1^4(k/d_{t-1})+2L_1^2(k/d_{t-1})$, where $L_t(\alpha)$ is the $(1-\alpha)$-th quantile of $N(0,(\zeta_tn)^{-1})$. These estimates are consistent according to \cite{sun2012scaled}. The algorithm is presented in Algorithm \ref{alg:scale}. As shown in Theorem \ref{thm:dt},  at time $t=\tau_k+\iota^*$, the problem reduces to an ordinary linear regression with high probability. 
Thus, during the hard selection stage, we estimate $\sigma_e$ by $\{C_y^{(t)}-2C_{Xy}^{(t)}(\widehat{J}_t)\widehat{\alpha}_t+\widehat{\alpha}_t^\top C_{X}^{(t)}(\widehat{J}_t)\widehat{\alpha}_t\}^{1/2}$.

\begin{onehalfspace}
\begin{breakablealgorithm}
  \renewcommand{\algorithmicrequire}{\textbf{Input:}}
  \renewcommand{\algorithmicensure}{\textbf{Output:}}
  \caption{Online Scaled Lasso}\label{alg:scale}
  \begin{algorithmic}[1]
  \STATE \textbf{For} $t=\iota^0, \iota^0+1,\ldots, $:
  \STATE \quad \textbf{input}: ${C}_y^{(t)},{C}_{Xy}^{(t)}(\widehat{J}_{t-1}),{C}_{X}^{(t)}(\widehat{J}_{t-1})$ and $\lambda_0$.
  \STATE \quad \textbf{initialize}: $\sigma_1\leftarrow\widehat{\sigma}_{\varepsilon,t-1} $, $\sigma_0\leftarrow0$, $\alpha_0\leftarrow\widehat{\alpha}_{t-1}$.
  \STATE \quad \textbf{while} $|\sigma_1-\sigma_0|>\delta_\sigma$:
  \STATE \quad \quad $\sigma_0\leftarrow\sigma_1$
  \STATE \quad \quad $\sigma_1\leftarrow (\zeta_tn)^{1/2} \left\|C_{y}^{(t)}-2C_{Xy}^{(t)}(\widehat{J}_{t-1})\alpha_0+\alpha_0^\top C_{X}^{(t)}(\widehat{J}_{t-1})\alpha_0\right\|_2$;
  \STATE \quad \quad $\lambda \leftarrow \sigma_1 \lambda_0$;
  \STATE \quad \quad $\alpha_0\leftarrow \min_\alpha L_\lambda\left({\alpha}\right)$;
  \STATE \quad \textbf{output}: $\widehat{\sigma}_{t} \leftarrow \sigma_1$.
  \end{algorithmic}  
\end{breakablealgorithm}
\end{onehalfspace}

\section{Theoretical Analysis}
\label{sec:theory}

In this section, we establish theoretical guarantees for the proposed RAVAS procedure in the
expanding-observable setting. Our analysis focuses on three key aspects:
(i) control of the dynamically selected model dimension,
(ii) coverage of the true active set by the selected model, and
(iii) estimation error of the final estimator once all important variables become observable.
To this end, we first introduce basic regularity conditions on the design and noise, and then
impose structural assumptions tailored to the presence of unobserved covariates and
time-varying pseudo-coefficients.

\subsection{Basic regularity conditions}

We begin with a restricted eigenvalue condition on the local sample covariance matrices
constructed at each change point.

\begin{assump}[Restricted eigenvalue condition]
\label{assump:RE}
For any $k\le K$, let
\[
J \subset \{p_{\tau_k-1}+1,\ldots, p_{\tau_k}\} \cup \Delta J,
\quad
\Delta J \subset \{1, \ldots, p_{\tau_k-1}\},
\quad
|\Delta J|\le \iota^0 n.
\]
Define
$C_{J,k}
=
(\iota^0 n)^{-1}\sum_{i=\tau_k}^{\tau_k+\iota^0 n} X_{i,J}^\top  X_{i,J}$.
We assume that, for any $k\le K$ and any such $J$, $C_{J,k}$ satisfies
the restricted eigenvalue condition $\mathrm{RE}(s_0)$, i.e., there exist
$s_0 = O(\log \Delta p)$, $\nu>0$, and $\phi>0$ such that for all
$J_0\subset J$ with $|J_0|\le s_0$ and all $\beta$ satisfying
$\|\beta_{J_0^c}\|_1\le \nu\|\beta_{J_0}\|_1$, it holds that
$\|\beta_{J_0}\|_1^2 \le \phi |J_0| \cdot \beta^\top C_{J,k} \beta$.
We refer to $s_0$ and $\phi$ as the restricted rank and eigenvalue of $C_{J,k}$, respectively.
\end{assump}

RE-type assumptions are known to be essentially unavoidable for reliable estimation and variable selection in high-dimensional regression.  
In the expanding-observable setting, Assumption~\ref{assump:RE} plays the same structural role.
Assumption~\ref{assump:RE} is required only at the $K$ time points $\tau_k+\iota^0$, rather than uniformly over all time points.  
Since $K$ is small relative to the data stream, the condition on this collection of sample covariance matrices is substantially less restrictive than global RE requirements.
What is more, Assumption~\ref{assump:RE} is imposed only on the matrices $C_{J,k}$, whose dimensions are moderate: each set $J$ contains at most $\Delta p$ newly observed variables together with at most $\iota^0 n = (\log\Delta p)^\kappa$ previously selected variables.  
Thus, the RE requirement concerns submatrices whose effective dimension is $O(\Delta p)$, rather than the full ambient dimension $p_t$.  
For such  subsets, RE-type conditions hold under broad covariance structures, including weak dependence, banded covariance, sparse-precision models and factor models, several of which will be illustrated in the
examples in Section \ref{subsec:nearly_sparse}.

In addition, since the online setting is most naturally modeled under a random design,
we impose the following distributional assumption on the covariates and noise.

\begin{assump}[Sub-Gaussian design and noise]
\label{assump:distribution}
The noise terms $\{\varepsilon_{ti}\}_{t,i}$ are i.i.d. sub-Gaussian random variables with mean zero and variance $\sigma_\varepsilon^2$.  
Similarly, the variables $\{Z_{ti}\}_{t,i}$ are i.i.d. continuous sub-Gaussian with mean zero and covariance $\Sigma_Z$.
\end{assump}


\subsection{Nearly sparsity and regression faithfulness}
\label{subsec:nearly_sparse}

When $t<\tau_K$, not all important covariates are observable. Following \citet{lu2012robustness} and \citet{lv2014model},
the misspecified sparse regression estimator satisfies
$E \big[ X_t (X_t^\top \alpha_*^{(t)} + U_t^\top \gamma_*^{(t)} - X_t^\top \alpha) \big] = 0$,
so the target pseudo-coefficient takes the form
\[
\alpha_{0}^{(t)}
=
\alpha_*^{(t)}
+
\big(\Sigma_{X}^{(t)}\big)^{-1}
\Sigma_{XU}^{(t)}
\,\gamma_*^{(t)},
\]
where $\Sigma_X^{(t)} = E[X_t X_t^\top]$, $\Sigma_U^{(t)} = E[U_t U_t^\top]$,
and $\Sigma_{XU}^{(t)} = E[X_t U_t^\top]$.
In the RAVAS procedure, estimation at time $t$ is carried out on the dynamically
selected submodel $\widehat{J}_{t-1}$. The relevant pseudo-coefficient becomes
\[
\tilde{\alpha}_{0}^{(t)}
=
\alpha_{*,\widehat{J}_{t-1}}^{(t)}
+
\{\Sigma_X(\widehat{J}_{t-1})\}^{-1}
\Sigma_{Z}(\widehat{J}_{t-1},U)\gamma_*^{(t)}.
\]

Two difficulties arise from this representation.
First, although $\alpha_*^{(t)}$ and $\gamma_*^{(t)}$ are sparse, the pseudo-bias term
$\{\Sigma_X(\widehat{J}_{t-1})\}^{-1}\Sigma_{Z}(\widehat{J}_{t-1},U)\gamma_*^{(t)}$
may be dense, challenging the feasibility of the regression.  
Thus, we require the {near-sparsity} of this pseudo-bias, which is ensured by a
separation structure on the covariance matrix.
Second, an important coordinate $j$ with $\alpha_{*,j}^{(t)}\neq 0$ may be masked by the
pseudo-bias, making $\tilde{\alpha}_{0,j}^{(t)}$ close to zero. Such cancellation
would cause false deletion, so a non-cancellation condition is needed. This is
analogous to regression faithfulness in graphical models.
To formalize these conditions, denote
\[
\overline{\eigenvalue}
=
\max_{|J|\le 2 \Delta p} \phi_{\max}(\Sigma_{Z,J}),
\qquad
\norm1 
=
\max_{|J|\le 2 \Delta p}\|\Sigma_{Z,J}\|_1,
\qquad
\underline{\eigenvalue}
=
\max_{|J|\le \iota^0 n} \phi_{\max}(\Sigma_{Z,J}),
\]
where $\phi_{\max}(A)$ is the maximal eigenvalue of $A$.

Drawing inspiration from the concept of faithfulness in causal discovery
graphical models, which implies that if $X$ and $Y$ have a direct causal relationship, they are not
conditionally independent given any subset of the graph, we first introduce the following regression faithfulness assumption to guarantee the non-cancellation of the signal.

\begin{assump}[Regression faithfulness]
\label{assump:faith}
Suppose that $\kappa$ is chosen such that
$(\log \Delta p)^{\frac{\kappa-1}{2}}
\gg
\underline{\eigenvalue}^{1/2}\norm1 \,\|\beta_*\|_2$.
For any $j\in J(\alpha_*^{(t)})$,
$|\tilde{\alpha}_{0,j}^{(t)}|
\gtrsim \faith$, 
where 
$\faith
\gg
\underline{\eigenvalue}^{1/2}\norm1 \,\|\beta_*\|_2
(\log \Delta p)^{-\frac{\kappa-1}{2}}$.
\end{assump}

Under growth conditions on $\|\beta_*\|_2$ and the eigenvalues of $\Sigma_Z$,
the lower bound $\faith$ can be chosen to tend to zero as $\Delta p\to\infty$.
Thus, Assumption~\ref{assump:faith} does not require signals to be uniformly bounded away from zero in dimension, but only rules out extreme cancellations between $\alpha_*^{(t)}$ and the pseudo-bias term.
Assumption~\ref{assump:faith} ensures that for any submodel containing $J(\alpha_*^{(t)})$, the signal strength in the presence of unobserved key variables is no smaller than order $\faith$. Similar conditions have been employed in existing works such as \citet{buhlmann2010variable} and \citet{ma2020variable}, where they are shown to be necessary for consistent model recovery. 
Although the condition appears technical, it excludes only a thin set of parameters:
consider
\[
\mathcal{M}_t=\left\{(\alpha^\top,\gamma^\top)^\top\in \mathbb{R}^{d_{t-1}+q_t}:
\exists j\in J(\alpha_*^{(t)}),\ 
\left|\left(\alpha+\{\Sigma_X(\widehat{J}_{t-1})\}^{-1}\Sigma_{Z}(\widehat{J}_{t-1},U)\gamma\right)_j\right|\ll\faith\right\},
\]
which is a thin slab in the $(d_{t-1}+q_t)$-dimensional space. Within any bounded parameter
region, the volume of $\mathcal{M}_t$ shrinks linearly to zero as $\faith\to 0$.
Thus, the requirement that the true parameter does not belong to $\mathcal{M}_t$ is generic, especially under the structured covariance settings considered in this subsection.
We next introduce the separation structure assumption to ensure that the pseudo coefficient is nearly sparse. 

\begin{assump}[Separation structure]
\label{assump:sigma}
Let $J_{\gamma_t} = J(\gamma_*^{(t)})$ and $\bar{\gamma}_t
= \gamma_{*,J_{\gamma_t}}^{(t)}$.  
For $J_t\subset\{p_{\tau_k-1}+1,\ldots, p_{\tau_k}\}\cup \Delta J$ with
$\Delta J\subset\{1, \ldots, p_{\tau_k-1}\}$, define 
$M_t
=
\Sigma_{J_t}^{-1} \Sigma_{J_t U}$.
There exists an index set $J_{str}\subset J_t$ with
$|J_{str}| = O(\log |J_t|)$ and a constant $0< c\le 1$ 
such that:
\begin{enumerate}[(i)]
\item \textbf{Strong rows:}
for any $j\in J_{str}$,
$\| M_t(j,J_{\gamma_t}) \|_2 \ge L_{str} \gtrsim  1$.
\item \textbf{Weak rows:}
for any $j\notin J_{str}$,
$\| M_t(j,J_{\gamma_t}) \|_1 \le L_{weak}
\ll \faith / \norm1$.

\item \textbf{Alignment or negligible projection:}
for any $j\in J_{str}$,
\[
|\langle M_t(j,J_{\gamma_t}),\bar{\gamma}_t\rangle|
\ge c \|M_t(j,J_{\gamma_t})\|_2 \|\bar{\gamma}_t\|_2
\quad \text{or} \quad
|\langle M_t(j,J_{\gamma_t}),\bar{\gamma}_t\rangle|
\ll \faith / \norm1.
\]
\end{enumerate}
\end{assump}

Conditions (i) and (ii) in Assumption~\ref{assump:sigma} require that the rows of $M_t$ can be classified into strong and weak categories according to their action on the support of $\gamma_*^{(t)}$. Consequently, the induced pseudo-signals can also be categorized into strong and weak types. The alignment condition in (iii) requires that, for each $j\in J_{str}$, the pseudo-signal $\langle M_t(j,J_{\gamma_t}),\bar\gamma_t\rangle$ is either of non-negligible order, proportional to $\|M_t(j,J_{\gamma_t})\|_2\|\bar\gamma_t\|_2$, or uniformly small. 
For most covariance structures, a row of $M_t$ is either substantially aligned with $\bar\gamma_t$, in which case the projection is automatically of the desired strong order; or nearly orthogonal to it, in which case the projection is small and the coordinate is effectively weak. 
The assumption simply excludes the case in which a strong row has a projection that lies in an intermediate range. Such intermediate cases require a delicate balance between $M_t(j,J_{\gamma_t})$ and $\bar\gamma_t$ and are structurally unstable in the sense that they disappear under small perturbations in the covariance or signal vector. 
Hence, (iii) holds robustly for a wide class of covariance models. We note that the assumption also permits the special case $L_{str}=0$, in which all rows of $M_t$ behave as weak rows. It corresponds to the stronger requirement that all relevant covariates are only weakly correlated with the unobserved variables. Then the pseudo-bias consists only of weak noise terms, making variable selection even more favorable.

We now present examples illustrating that
Assumption~\ref{assump:sigma}(i)–(ii) hold under widely used covariance structures.

\begin{eg}[Banded covariance]
\label{eg:banded}
Suppose $\Sigma_Z$ is $M$-banded, i.e.,
$(\Sigma_Z)_{ij}=0$ whenever $|i-j|>M$.
By the classical result of \citet{demko1984decay},
the inverse $\Sigma_{J_t}^{-1}$ exhibits exponential
off-diagonal decay, so in particular
$\|\Sigma_{J_t}^{-1}\|_1$ is bounded uniformly in $t$.
Likewise, $\Sigma_{J_t U}$ is banded with the same bandwidth.
Hence each row of
$M_t = \Sigma_{J_t}^{-1}\Sigma_{J_t U}$
concentrates its mass near the coordinates of $J_{\gamma_t}$.
Rows whose indices lie within distance $M$ of $J_{\gamma_t}$
form $J_{str}$ and satisfy
Assumption~\ref{assump:sigma}(i), while distant rows
satisfy the weak bound in (ii).
\end{eg}

\begin{eg}[Factor model]
If $Z$ is generated by the following factor model,
$Z=\Lambda f + \epsilon$,
where $\Lambda\in\mathbb{R}^{p_Z\times r}$ is the loading matrix, the factor $f$
satisfies $E[ff^\top]=I_r$, and $\epsilon$ is $p_Z$-dimensional noise with
mean zero and $\cov(\epsilon)=\Psi$, where $\Psi$ is diagonal with finite entries.
Then
$M_t
=
\Psi_{J_t}^{-1}\Lambda_{J_t}Q_t\Lambda_U^\top,
\quad
Q_t
=
\big(I_r+\Lambda_{J_t}^\top\Psi_{J_t}^{-1}\Lambda_{J_t}\big)^{-1}$.
If the row norms of $\Lambda$ differ substantially (e.g., strong-loading rows and
weak-loading rows), then the corresponding rows of $M_t$ inherit this strong–weak
structure. Strong-loading rows satisfy Assumption~\ref{assump:sigma}(i), while
weak-loading rows produce small $\ell_1$ norms, satisfying (ii).
\end{eg}

\begin{eg}[Sparse precision]
Assume $Z \sim N(0, \Sigma_{Z})$ and $\Theta_{Z} = \Sigma_{Z}^{-1}$ corresponds
to a sparse undirected graph $G$ with bounded degree and walk-summability, so
$|(\Sigma_{Z})_{ij}| \leq C \rho^{\mathrm{dist}_{G}(i,j)}, \ 0 < \rho < 1$,
and hence $\|\Sigma_{Z}\|_{1} \leq C' < \infty$.
In addition, suppose there exists a separator set $S \subset \{1, \ldots, p_{Z}\}$ with $|S| = O(\log \Delta p)$ such that, for all $t$, any path from $X_{t}\setminus S$ to $U_{t}$ in $G$ must pass through $S$. As in the short-range case, walk-summability implies that the entries of
$M_{t} = \Sigma_{\widehat{J}_{t-1}}^{-1} \Sigma_{\widehat{J}_{t-1} U}$
decay with graph distance, and it can be verified that $M_{t}$ satisfies Assumption~\ref{assump:sigma}(i) and (ii).
\end{eg}

These examples demonstrate that the separation structure arises naturally under common
covariance models, ensuring that the pseudo-bias term $M_t\gamma_*^{(t)}$ is
near-sparse and enabling reliable variable selection under expanding observability. It is worth noting that continuously varying covariance structures are not suitable in this context, such as those with smooth exponential decay. Under such models, the pseudo-bias term tends to diffuse its mass across many coordinates, failing to produce a clear contrast between (fake) signals and weak noise. As a result, the pseudo-coefficient is no longer approximately sparse, and lasso-based variable selection cannot reliably separate important variables from the accumulation of many small but non-negligible effects. The separation structure is therefore essential for enabling effective screening under expanding observability.

\subsection{Main results: model size, coverage, and estimation error}
\label{subsec:main_results}

The core idea of RAVAS is to perform phased and batched dimension reduction on the model, keeping the selected model size at a manageable scale while preserving the important variables. This, in turn, allows us to control both the coverage probability and the estimation error in an online, expanding-observable environment. 
We first impose a sparsity and minimal signal condition on the underlying true parameter.

\begin{assump}
\label{assump:signal}
The true active set has size $s_* = |J(\beta_*)| = O(\log \Delta p)$.  
There exists a constant $\bar{C}_\beta<\infty$ such that $\|\beta_*\|_\infty\le \bar{C}_\beta$, and the minimal signal strength satisfies
$\min_{j\in J(\beta_*)} |\beta_{*,j}| \gtrsim \faith$.
\end{assump}

Under Assumption~\ref{assump:sigma}, for any $t$, the pseudo-coefficient $\tilde{\alpha}_0^{(t)}$ is nearly sparse with threshold of order $\faith$. Define
$H_{\faith}(\tilde{\alpha}_0^{(t)})
=
\{j\in\widehat{J}_{t-1}: |\tilde{\alpha}_{0,j}^{(t)}|\gtrsim \faith\}$
 and 
$s_{0,t} = |H_{\faith}(\tilde{\alpha}_0^{(t)})|$.
The quantity $s_{0,t}$ represents the effective sparsity of the pseudo-coefficient at time $t$, allowing for additional ``fake'' signals induced by temporarily unobserved variables. Our first result shows that RAVAS controls the selected model dimension around $s_{0,t}$ in each cycle.

\begin{thm}
\label{thm:dt}
Suppose that Assumptions \ref{assump:RE}--\ref{assump:signal} hold and $\lambda_t^0$ is selected by Algorithm~\ref{alg:RAVAS}. For any $\tau_k\le t<\tau_{k+1}$, if $\tau_{k+1}-\tau_k\le \iota^*$, then
\[
P\big\{d_{\tau_{k+1}-1}=O\big(\underline{\eigenvalue}\,s_{0,t}\big)\big\}
=
1-O\big(k(\iota^0 n)^{-c_1}\big).
\]
Otherwise,
\[
P\big(d_{\tau_{k+1}-1}=s_{0,t}\big)
=
1-O\big(k\Delta p^{-c_2}\big).
\]
\end{thm}

Theorem~\ref{thm:dt} shows that, in each cycle, the model dimension $d_t$ is effectively controlled by the effective sparsity $s_{0,t}$. When the cycle is too short to trigger hard selection ($\tau_{k+1}-\tau_k\le \iota^*$), soft selection alone reduces the model size to order $O(\underline{\eigenvalue}\,s_{0,t})$ with high probability, where $\underline{\eigenvalue}$ reflects the local conditioning of the design. When the cycle is long enough and the hard selection step is executed, the model is further refined and $d_t$ matches $s_{0,t}$ up to constants. 

Regarding $s_{0,t}$, note that when some important variables are still unobserved ($k<K$), the pseudo-coefficient $\tilde{\alpha}_0^{(t)}$ may contain fake signals arising from the bias term $M_t\gamma_*^{(t)}$. Under the separation structure, these fake signals are few in number, i.e. $s_{0,t}=O(\log \Delta p)$ when $k<K$.
Once all important variables have become observable ($k\ge K$), the pseudo-coefficient aligns with the true coefficient, i.e. $s_{0,t}=s_*$. As a consequence, RAVAS maintains a model dimension of order $O(\log \Delta p)$ before full observability and of order $s_*$ afterwards.

Our next result establishes a coverage property for the selected model.

\begin{thm}
\label{thm:Pcover}
Suppose that Assumptions \ref{assump:RE}--\ref{assump:signal} hold and $\lambda_t^0$ is selected by Algorithm~\ref{alg:RAVAS}. Then, for any $\tau_k\le t<\tau_{k+1}$,
\[
P\big(J(\alpha_*^{(t)})\subset\widehat{J}_{t}\big)
=
1-O\big(k\Delta p^{-c_3}\big).
\]
\end{thm}

The main idea of Theorem~\ref{thm:Pcover} is that, with high probability, the dynamically selected model contains the support of the  pseudo-coefficient at every step of the procedure. Combined with the regression faithfulness condition, this ensures that all truly important coordinates are retained and are never removed by the online screening mechanism. This guarantee holds throughout the entire algorithmic process, during both soft- and hard-selection phases and regardless of where each cycle terminates. In particular, after $\tau_K$, when $\alpha_*^{(t)}$ coincides with $\beta_*$, the selected model continues to cover the true active set $J(\beta_*)$ with high probability.

We finally turn to the estimation error of the RAVAS estimator. Since, for $t<\tau_K$, the target is a pseudo-coefficient affected by missing variables, the resulting error bound is less informative about the underlying mechanism. We therefore focus on the regime $t\ge \tau_K$, where all important covariates are observable and the target is the true coefficient $\beta_*$.

\begin{thm}
\label{thm:err}
Suppose that Assumptions \ref{assump:RE}--\ref{assump:signal} hold and $\lambda_t^*$ is selected as in \eqref{eq:lam opt raoa}. Let $k \ge K$, $t\ge \tau_k$ and $\zeta_t$ is defined in Section~\ref{sec:model}.
\begin{itemize}
\item[(1)] When $\tau_{k}+\iota^0 \le t <\tau_{k}+\iota^*$ (soft-selection stage),
\[
\|\widehat{\alpha}_{t}^*-\beta_{*}\|_1
=
O\left(s_*\sqrt{\frac{\log d_{t-1}}{\zeta_tn}}\right)
\]
holds with probability $1-O\big\{d_{t-1}^{-1}+(k-1)\Delta p^{-c_4}\big)$.
\item[(2)] When $\tau_{k}+\iota^*\le t < \tau_{k+1}$ (hard-selection stage),
\[
\|\widehat\alpha_t^*-\beta_{*}\|_2^2
=
O\big\{s_* (\zeta_tn)^{-1}\big\}
\]
holds with probability $1-O\big(k\Delta p^{-c_5}\big)$.
\end{itemize}
\end{thm}

Theorem~\ref{thm:err} shows that, once all important variables have become observable, RAVAS operates within a substantially reduced model space with $d_{t-1}\ll p_t$, and hence yields estimation errors far smaller than those attainable in classical high-dimensional regression. During the soft-selection period, the $\ell_1$ error scales as 
$O\{s_* \sqrt{(\log d_{t-1})/(\zeta_tn)}\}$, and after hard selection reduces the model to order $s_*$, the estimator achieves the optimal $\ell_2$ rate $O(s_* (\zeta_tn)^{-1})$, matching the performance of an oracle that knows the true support.
These guarantees come with slightly weaker probability bounds, a reflection of the online and dynamically evolving nature of the feature space. However, this trade-off is acceptable: the procedure attains substantially improved accuracy while avoiding computations at the full ambient dimension, which is infeasible in the expanding-observable setting. Indeed, such mitigation is unavoidable, as processing all variables at every time step would exceed both memory and computational budgets. 

We close this section by noting that the theoretical results also cover the degenerate
fully observable setting, for which several arguments simplify due to the absence of
unobserved covariates. Related discussions are provided in Section 3 of the Supplementary Materials.

\section{Numerical Experiments}
\label{sec:simu}

We simulate the expanding-observable case using 100 replicates. Each data block contains $n=50$ observations, and the total number of data blocks is set to $T_{max}=4000$, resulting in a full sample size of $N_{max}=2\times 10^5$. New variables emerge at the following time points,
$$
\tau_k=\left\{
\begin{aligned}
&50\cdot k,\ \ k\le 20;\\ 
&1000+(k-20)\times 100, 20<k\le 24; 
\end{aligned}
\right.
$$
and additionally $\tau_{25},\tau_{26},\tau_{27}=8\times1600,2000,2800$. At each $\tau_k$, $\Delta p=500$ newly-emerged features are introduced, starting from an initial model dimension of $p_1=500$. By the end, there are $p_{T_{\max}}=14000$ features. Let $N_t$ represent the observed sample size up to time $t$. When $t\le1000$, the ratio $p_t/N_t=0.2$, reflecting an expanding-observable setting. When $t>1000$, the gaps between the change points become larger, indicating that the system tends to stabilize. 
The underlying full model consists of $s=10$ important variables, with values randomly selected from $\{-10,\ldots,-1,1,\ldots,10\}$. At time $t=1$, five important features are observed, resulting in $J(\alpha_*^{(1)})=5$. Without loss of generality, we assume $\alpha_{*,i}^{(1)}\neq0$ for $i\le 5$. At times $\{50,100,250,600,1600\}$, one more important feature is observed, specifically the first among the newly-emerged variables. The covariate $X_t$ are $m_t$ blockwise dependent, where $m_t=0.02p_t$ increases with  $p_t$.  Specifically, {variables sharing the same remainder when divided by 50 are correlated}. The covariance structure for the correlated features is exchangeable, with non-diagonal entries set to $\rho=0,\ 0.3\ ,0.5$ and $0.7$. The noise level is $\sigma_\varepsilon=1$. {In this simulation, we evaluate a range of values $\{3,6,10\}$ for $\iota^0$. The results, presented in the Supplementary Material, demonstrate that the RAVAS method is insensitive to the choice of $\iota^0$ under banded covariance structure. For illustration in the main text, we select $\iota^0=6$ to showcase the performance of RAVAS. The other parameter is configured as $\iota^*=\Delta p/n=20$.}

\begin{figure}[htp]
\centering
\includegraphics[width=0.7\textwidth]{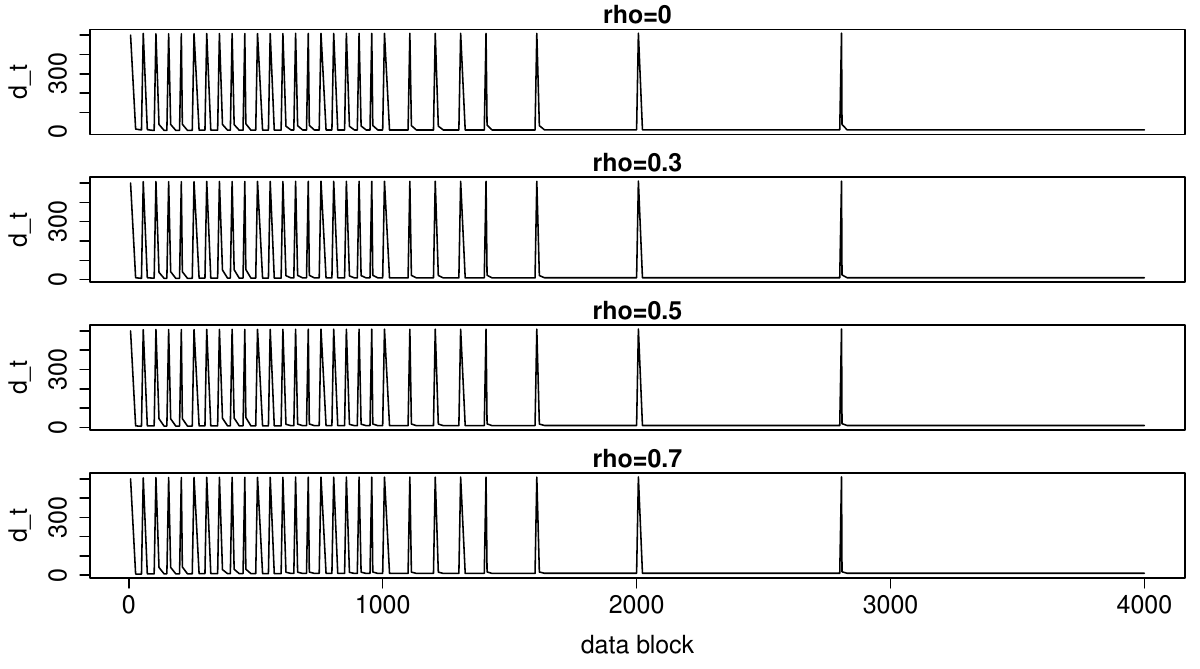}
\caption{\label{fig:dt id} The  numbers of selected variables of the proposed RAVAS method versus the number of observed data blocks. }
\end{figure}

\begin{figure}[htp]
\centering
\includegraphics[width=0.7\textwidth]{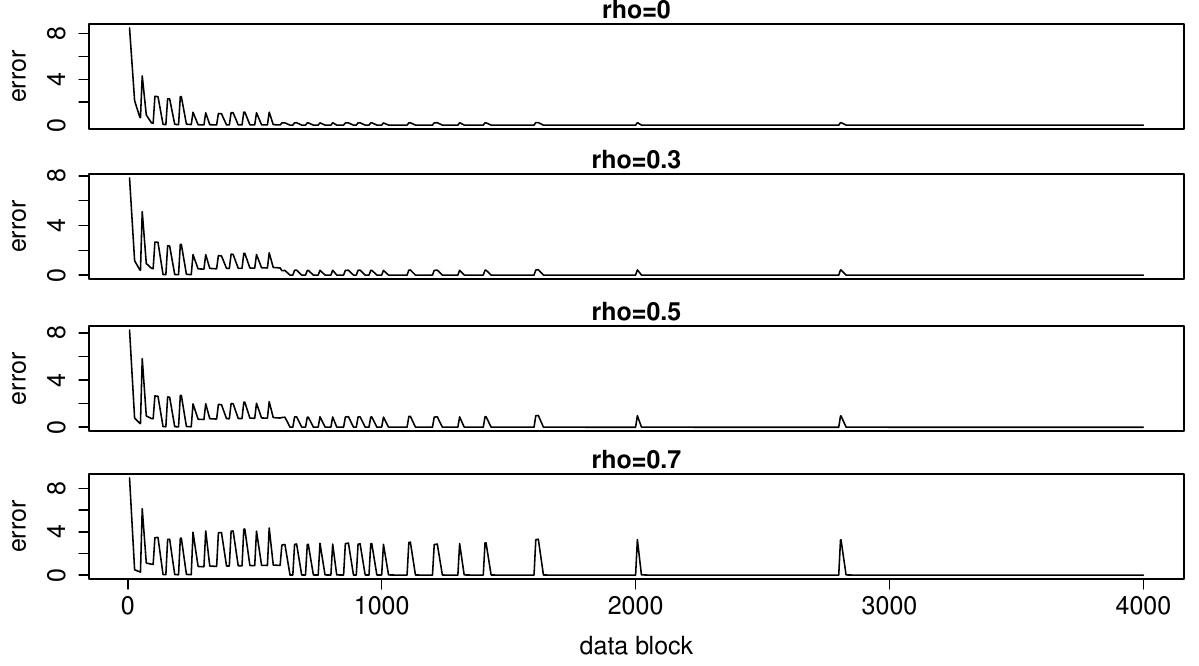}
\caption{\label{fig:err id} The estimation error of the proposed RAVAS method versus the number of observed data blocks. }
\end{figure}

\begin{figure}[htp]
\centering
\includegraphics[width=0.7\textwidth]{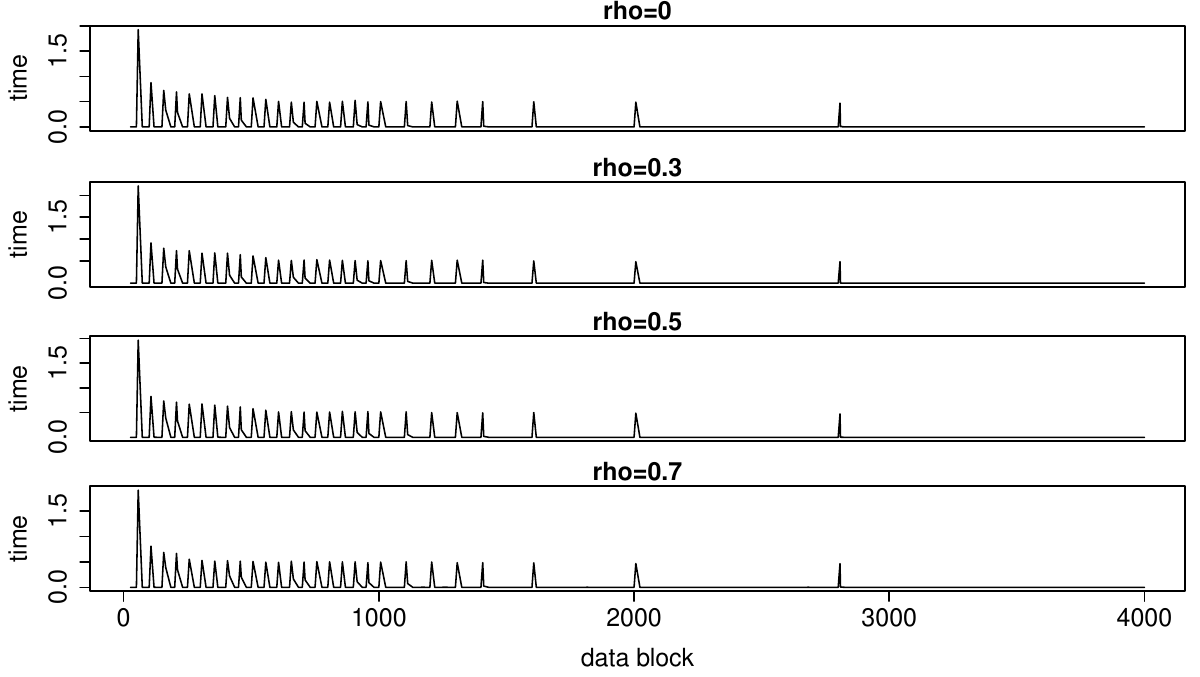}
\caption{\label{fig:time id} The computing time for each update  of the proposed RAVAS method versus the number of observed data blocks.}
\end{figure}

As shown in Figure \ref{fig:dt id}, the number of selected variables by the proposed method changes according to the emergence of new variables, which is desired in the expanding-observable setting. At the start of each period, the dimension is of order $O(\Delta p)$, but it rapidly decreases to $O(s_t)$ as more data are observed.  Figure \ref{fig:err id} indicates that the estimation error $E\|\widehat{\alpha}_t-\beta_{*,\widehat{J}_t}\|_2^2$ is dynamically changing along with the changes of the system.  In each panel, the peak value gradually decreases, reflecting a reduction in the variance of pseudo noise $\sigma_e^2$. The final estimation errors at time $T_{max}$ for $\rho=0,\ 0.3,\ 0.5,\ 0.7$ is $1.6\times 10^{-4},\ 2.2\times 10^{-4},\ 3.0\times 10^{-4},\ 5.1\times 10^{-4}$, respectively. The computing time for each step follows a similar pattern to that of the model dimension, as illustrated in Figure \ref{fig:time id}.

In addition to evaluating performance under expanding observability, we also consider the fully observable setting and compare RAVAS with existing online variable selection methods designed for static feature spaces. Detailed results and discussions can be found in Section 3 of the Supplementary Materials.

\section{Data Applications} 
\label{sec:realdata}

\subsection{Analysis of factors influencing PM2.5}

In this example, we utilize the RAVAS method to explore the factors influencing daily PM2.5 concentrations in Beijing. The response variable is derived from PM2.5 concentrations based on near-surface visibility recorded at 4,011 sites across the Northern Hemisphere \citep{pm25}, accessible at \url{https://data.tpdc.ac.cn/zh-hans/data/d3878788-7bef-4249-9973-07d4b2e46cb8}. We calculate the average daily PM2.5 levels at sites located in Beijing (latitude $37.0^\circ-43.2^\circ$N, longitude $110^\circ-120^\circ$E) over the period from 1991 to 2018. The predictors are sourced from two datasets, both available for download at \url{https://data.tpdc.ac.cn/home}. The first is the China Meteorological Forcing Dataset (CMFD), which includes seven near-surface meteorological elements \citep{cmfd}: daily 2-meter air temperature, surface pressure, specific humidity, 10-meter wind speed, downward shortwave radiation, downward longwave radiation, and precipitation rate, all spanning the years 1991 to 2018. The second dataset, the High-Resolution and High-Quality Pollutants Dataset for China (CHAP), contains concentrations of various pollutants recorded from different starting years \citep{chap}: O3 and PM10 from 2000 to 2018; NO2 from 2008 to 2018; and CO and SO2 from 2013 to 2018, which demonstrates a typical expanding-observable scenario. The meteorological factors -- temperature, humidity, pressure, and wind speed -- are often weakly correlated or nearly uncorrelated, reflecting a sparse dependency structure. In contrast, pollutant concentrations tend to exhibit localized correlations; for example, NO2 and SO2 are often strongly correlated due to shared emission sources, whereas their correlation with O3 is typically weaker owing to differing chemical processes \citep{cui2018relationship,ling2019correlation}. These meteorological factors and pollutant concentrations are used to predict PM2.5 levels in Beijing. 

For each year, we utilize meteorological factors from the previous day and pollutant concentrations from the same day to predict PM2.5 levels, resulting in 364 observations per year. We use data from 1991 to 2017 as the training set (sample size $N=9828$) and the data from 2018 as the test set (sample size $M=364$). We treat weekly data as a data block, setting $n=7$, which results in a total of 1,404 data blocks. Our analysis focuses on whether meteorological elements and pollutant concentrations from adjacent regions influence PM2.5 levels in Beijing. We divide the latitude into 15 equal segments and the longitude into 20 equal segments, yielding a total of 300 grids. Due to the lack of monitoring sites, some features contain missing values in certain regions. The model dimensions are $p_t=2016, 2594, 2882, 3458$ when $t=1,417,781,990$ (the first week of year 1991, 2000, 2008, 2013). We set {$\iota^0=100$ and $\iota^*=\lceil 2\Delta p/n \rceil$}, where $\lceil\cdot\rceil$ denotes rounding up to the nearest integer.

\begin{figure}[htp]
\centering
\includegraphics[width=0.6\textwidth]{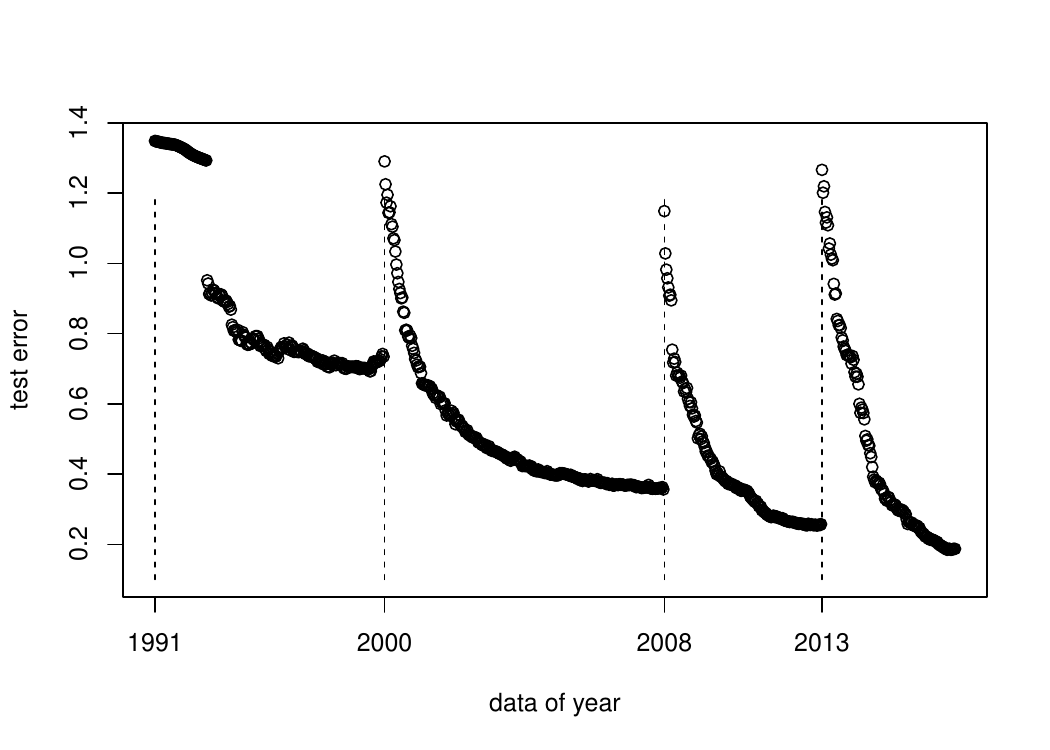}
\vspace{-0.2in}
\caption{\label{fig:err pm2.5} Test errors over observation time in the PM2.5 dataset, which shows a periodic pattern. }
\end{figure}

\begin{figure}[htp]
\centering
\includegraphics[width=1\textwidth]{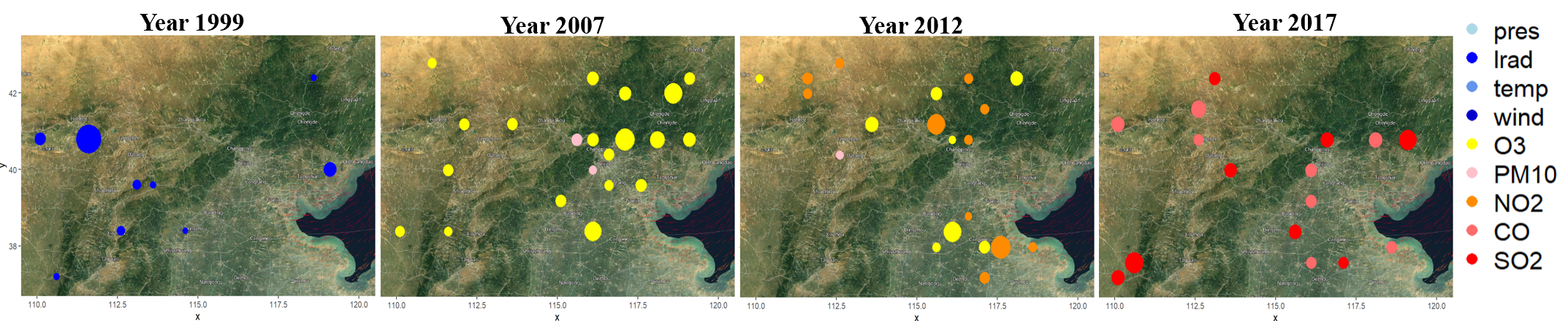}
\vspace{-0.2in}
\caption{\label{fig:map}  Selected important features at the end of the years 1999, 2007, 2012, and 2017. The colors of the dots represent different pollutants, while the dot sizes correspond to the magnitude of the coefficient values. }
\end{figure}

The test errors versus data blocks are plotted Figure \ref{fig:err pm2.5}. The testing errors can be categorized into four stages corresponding to changes in the system, with a consistent decreasing trend in each stage. As the number of variables increases, the errors also decrease incrementally. Based on the data up to 2017, the error reaches a minimum of 0.184. We further illustrate the important features according to their geographical locations at the end of each stage in Figure \ref{fig:map}. {Notably, meteorological factors predominantly exhibit a negative influence, while pollutant concentrations have a primarily positive impact, corresponding to cool tones (blue) and warm tones (yellow, red, etc.), respectively. The relationships among NO2, CO, SO2, and PM2.5 have been extensively studied in the literature \citep{li2019spatiotemporal, xie2015spatiotemporal}, aligning with our final findings. The RAVAS method reveals two notable phenomena. First, as pollutant concentrations are incorporated into the model, the influence of meteorological factors gradually weakens. Second, following the introduction of NO2 concentrations recorded from 2008, the significance of O2 concentrations diminishes. This pattern reflects the exploratory nature of scientific discovery, where knowledge evolves alongside technological advancements and more comprehensive data collection. Furthermore, this process provides a possible validation mechanism: when the addition of new features no longer alters the model structure, it suggests that the model may have reached its correct specification. Such insights cannot be attained in a fully observable setting, highlighting the advantages of the expanding observable framework.} 

\vspace{-0.2in}

\subsection{Peer-to-peer credit evaluation dataset}

The second data example is the peer-to-peer credit evaluation dataset from the leading European lending platform, Bondora, which can be downloaded from \url{https://www.kaggle.com/datasets/sid321axn/bondora-peer-to-peer-lending-loan-data/data}. 
This file consists of loan records of 1,79,236 borrowers. The time period of the loan is from Feb 2009 to July 2021. 
Although our theoretical analysis focuses on continuous covariates, this real-data experiment illustrates that the proposed method performs robustly even when the feature space includes mixed continuous and discrete variables.
We preprocess the Bondora peer-to-peer lending data by first removing identifiers, textual fields, and all post-origination or outcome-related variables. 
The raw data are preprocessed by standardizing continuous variables, encoding categorical attributes in a compact form, and simplifying high-cardinality fields (e.g., job titles) into grouped categories. Only pre-origination borrower and loan characteristics are retained to ensure validity for predictive modeling.

We order the Bondora dataset by borrower ID and treat it as a streaming sequence. The data contains 2261 blocks in total. Each block contains 922 obsevations. We use first 2250 as training and the last 11 blocks as the testing data. We set $\iota^0n=922$ and $c_h=1$.
At the beginning of the stream, each observation contains a core set of borrower, credit-history, and loan-application variables. As the stream progresses, additional credit-bureau and account-level indicators gradually become observable. 
The expanding availability of these variables leads to several feature-increment events, each introducing a block of new predictors. When combined with selected interaction terms, the effective model dimension increases from roughly 3,000 variables in the initial stage to nearly 8,000 variables after all increments. This creates a realistic high-dimensional online setting in which the covariate space grows substantially over time, providing a natural environment of expanding observability.

As shown in Table \ref{tab:P2P}, the model dimension, computing time, and prediction error all decrease as data collection progresses.  
Across the stages, the selected features exhibit a clear progression from broad to focused structure. Initially, the model relies on generic credit-risk indicators and coarse borrower and loan characteristics, along with the interaction terms. As more information becomes available, the selection gradually shifts toward more refined features that capture subtle differences across regions and borrower groups. In the final stage, the model stabilizes on a relatively compact set of core risk measures combined with a small number of targeted interaction effects, reflecting a parsimonious yet expressive representation of the interest-rate mechanism.

\begin{table}[htbp]
  \centering
  \caption{\label{tab:P2P} Numbers of selected variables, computational time and prediction errors ($Err_t=M^{-1}\sum_{i=1}^M(y_i-\widehat{y}_i)^2$) at different observation time points (data blocks), where $M$ is the sample size of testing data.}
    {\small\begin{tabular}{|l|c|c|c|c|c|c|c|c|}
    \hline
    \textbf{data block }  & 100 & 500 & 800 & 1300 & 2000  & 2250 \\
    \hline
    \textbf{No. selected variables}   &  361 & 65 & 64 & 65  & 36  &  25 \\
    \hline
    \textbf{computational time} & 0.2034 &0.0055& 0.0055 &0.0056& 0.0028 &0.0019 \\
    \hline
    \textbf{prediction error} & 0.2987 & 0.1545& 0.1110 &0.0097 &0.0009 &0.0005\\
    \hline
    \end{tabular}}%
\end{table}%

\section{Discussion}

In this work, we develop the RAVAS method, addressing the significant challenges of high-dimensional regression in online settings with increasing feature dimensions. Our approach combines soft and hard selection techniques to effectively identify important predictors while managing computational efficiency. We establish theoretical guarantees that ensure the capturing of crucial features with high probability. The innovations introduced not only offer practical solutions for real-time high-dimension data analysis, but also pave the way for future research in adaptive statistical methodologies for expanding-observable problems.

Future investigations could explore the extension of the RAVAS framework to scenarios involving dynamic underlying coefficients, where the relationships between predictors and the response may change over time. This would necessitate the development of adaptive methods that not only update model parameters but also account for temporal fluctuations in the significance of predictors. Additionally, incorporating more sophisticated methods for handling unobserved important data such as instrumental variable method could further enhance the accuracy and robustness  of the model.

\bigskip
\begin{center}
{\large\bf SUPPLEMENTARY MATERIAL}
\end{center}

The supplementary materials provide technical details, algorithmic components, and additional numerical evidence that support the results presented in the main paper.

\begin{center}
{\large\bf Acknowledgement}
\end{center}
This research is supported by the National Natural Science Foundation of China (No. 12292981, 12571309, 62588101, 12288101), National Key Research and Development Program
of China grants (No. 2022YFA1003800), the New Cornerstone Science Foundation through Xplorer Prize, the Fundamental Research Funds for the Central Universities (LMEQF), an EPSRC grant EP/W014971/1, and the LMAM. 

\begin{center}
{\large\bf Data Availability}
\end{center}
The code can be found on our GitHub page at \url{https://github.com/anneyang0060/ExpandingObservables}. We have included detailed instructions in the ‘Readme’ file on how to reproduce the numerical results.

\begin{center}
{\large\bf Disclosure Statement}
\end{center}
The authors report there are no competing interests to declare.

{ 
\bibliographystyle{agsm}
\bibliography{ref}
}
\end{document}